\documentclass[final]{birkmult}

\numberwithin{equation}{section}

\theoremstyle{plain}
\newtheorem{theorem}[equation]{Theorem}
\newtheorem{lemma}[equation]{Lemma}
\newtheorem{corollary}[equation]{Corollary}

\theoremstyle{definition}
\newtheorem{definition}[equation]{Definition}
\newtheorem{assumption}[equation]{Assumption}
\theoremstyle{remark}

\usepackage{amsfonts,amssymb,latexsym}  
\usepackage{amsmath,enumerate}

\usepackage[authoryear]{natbib}

\numberwithin{equation}{section}

\newcommand{\N}{\ensuremath{\mathbb{N}}}    
     
\newcommand{\Z}{\ensuremath{\mathbb{Z}}}    
\renewcommand{\P}{\ensuremath{\mathbb{P}}}

\newcommand{\E}{\ensuremath{\mathbb{E}}}

\newcommand{\Id}{\mathbf {1}}

\newcommand{\TT}{\mathcal T}

\newcommand{\ee}{e^{\frac1e}}
\newcommand{\given}{\ \vert \ }

\newcommand{\suf}{{\text{suf}}}
\renewcommand{\tau}{\TT}

\begin{document}
%
%
%
%
%
%
%
%
%
\title[Inequalities for Empirical Context Trees]
 {Exponential inequalities for empirical unbounded context
  trees}
\author[A. Galves]{Antonio Galves}

\address{%
  Instituto de Matem\'atica e Estat\'{\i}stica\\
   Universidade de S\~ao Paulo\\
   BP 66281, 05315-970\\
   S\~ao Paulo, Brasil}

\email{galves@ime.usp.br}

\thanks{This work is part of
PRONEX/FAPESP's project \emph{Stochastic behavior, critical phenomena
  and rhythmic pattern identification in natural languages} (grant
number 03/09930-9) and CNPq's projects \emph{Stochastic modeling of speech}
(grant number 475177/2004-5) and  \emph{Rhythmic patterns, prosodic domains and probabilistic modeling in Portuguese Corpora} (grant number 485999/2007-2). AG is partially supported by a CNPq fellowship (grant 308656/2005-9)
and FL is supported by a FAPESP fellowship (grant 06/56980-0)}
\author[F. Leonardi]{Florencia Leonardi}
\address{
 Instituto de Matem\'atica e Estat\'{\i}stica\\
   Universidade de S\~ao Paulo\\
   BP 66281, 05315-970\\
   S\~ao Paulo, Brasil}
\email{leonardi@ime.usp.br}



\date{January 1, 2004}

\begin{abstract}
In this paper we obtain non-uniform exponential upper bounds for the rate of convergence of a version of the 
algorithm Context, when the underlying tree is not necessarily bounded. 
The algorithm Context is a well-known tool to estimate the context tree of a Variable Length Markov Chain. 
As a consequence of the exponential bounds we obtain a strong consistency result. We generalize in this way several previous results in the field.   
\end{abstract}

\maketitle

\section{Introduction}

In this paper we present an exponential bound for the rate of convergence 
of the algorithm Context for a class of unbounded variable memory models, taking values on a finite 
alphabet $A$. From this it follows a strong consistency result for the algorithm Context in this setting. 
Variable memory models were first introduced in the information 
theory literature by Rissanen \cite{rissanen1983} as a universal system for data compression.
Originally  called by Rissanen {\sl   finite memory source} or {\sl probabilistic tree},
this class of models recently became popular in the statistics 
literature under the name of {\sl Variable Length Markov Chains   (VLMC)} 
\cite{buhlmann1999}.

The idea behind the notion of variable memory models is that the probabilistic
definition of each symbol only depends on a finite part of the past and the length 
of this relevant portion is a function of the past itself. Following Rissanen we 
call \emph{context} the minimal relevant part of each past. The set of all contexts satisfies the 
suffix property which means that no context is a proper suffix of another context. This property allows to represent the set of all contexts as a rooted labeled tree. With this representation the process is described by the tree of all contexts and a associated family of probability measures on $A$, indexed by the tree of contexts. Given a context, its associated probability measure gives the probability of the next symbol for any past having this context as a suffix. From now on the pair composed by the context tree and the associated family of probability measures will be called \emph{probabilistic context tree}.

Rissanen not only introduced the notion of variable memory models 
but he also introduced the algorithm Context to estimate the probabilistic
context tree. The way the algorithm Context works can be summarized as follows. Given 
a sample produced by a chain with variable memory,  we start with a maximal tree of  candidate 
contexts for the sample. The branches of this first tree are then 
pruned until we obtain a minimal tree of contexts well adapted to the 
sample. We associate to each context  an 
estimated probability transition defined as the proportion of time the context appears in the sample followed by each one of the symbols in the alphabet. From Rissanen \cite{rissanen1983} to Galves et al. \cite{galves2006}, passing by 
Ron et al. \cite{ron1996} and B\"uhlmann and Wyner \cite{buhlmann1999}, several variants of the 
algorithm Context have been presented in the literature. In all the variants the decision to prune a branch is taken by considering a {\em cost} function. A branch is pruned if the cost function assumes a value 
smaller than a given threshold. The estimated context tree is the 
smallest tree satisfying this condition. The estimated family of probability 
transitions is the one associated to the minimal tree of contexts.

In his seminal paper Rissanen proved the weak consistency of the algorithm Context 
in the case where  the contexts have a bounded length, i.e. where the tree of contexts is finite. B\"uhlmann and Wyner \cite{buhlmann1999} proved the weak consistency of the algorithm also in the finite case without assuming a priori known bound on the maximal length of the memory, but using a bound allowed to grow with the size of the sample. In both papers the cost function is defined using the log likelihood ratio test to compare two candidate trees and the main ingredient of the consistency proofs was the chi-square approximation to the log likelihood ratio test for Markov chains of fixed order. A different way to prove the consistency in the finite case was introduced in \cite{galves2006}, using exponential inequalities for the estimated transition probabilities associated to the candidate contexts. As a consequence they obtain an exponential upper bound for the rate of convergence of their variant of the algorithm Context.

 The unbounded case as far as we know was first considered by Ferrari and Wyner \cite{ferrari2003} who also proved a weak consistency result for the algorithm Context in this more general setting. The unbounded case was also considered by Csisz\'ar and Talata \cite{csiszar2006} who introduced a different approach for the estimation of the probabilistic context tree using the Bayesian Information Criterion (BIC) as well as the Minimum Description Length Principle (MDL). We refer the reader to this last paper for a nice description of other approaches and results in this field, including the context tree maximizing algorithm by  Willems et al. \cite{willems1995}. With exception of Weinberger et al. \cite{weinberger1995}, the issue of the rate of convergence of the algorithm estimating the probabilistic context tree was not addressed in the literature until recently.  Weinberger et al. proved in the bounded case that the probability that the estimated tree differs from the finite context tree generating the sample is summable as a function of the sample size. Duarte et al. in \cite{duarte2006} extends the original weak consistency result by Rissanen \cite{rissanen1983} to the unbounded case. Assuming weaker hypothesis than \cite{ferrari2003}, they showed that the on-line estimation of the context function decreases as the inverse of the sample size.

In the present paper we generalize the exponential inequality approach presented in \cite{galves2006} to obtain an exponential upper bound for the algorithm Context in the case of unbounded probabilistic context trees. Under suitable conditions, we prove that the truncated estimated  context tree converges exponentially fast to the tree generating the sample, truncated at the same level. This improves all results known until now.

The paper is organized as follows. In section~\ref{def} we give the 
definitions and state the main results. Section~\ref{proof1} is devoted to the 
proof of an exponential bound for conditional probabilities, for unbounded probabilistic context trees. In section \ref{proofmain}  we apply this exponential bound to estimate the rate of 
convergence of our version of the algorithm Context and to prove its consistency.

\section{Definitions and results}\label{def}

In what follows $A$ will represent a finite alphabet of size $|A|$.
Given two integers $m\leq n$, we will denote by $w_m^n$ the sequence
$(w_m, \ldots, w_n)$ of symbols in $A$. The length of the sequence
$w_m^n$ is denoted by $\ell(w_m^n)$ and is defined by $\ell(w_m^n) =
n-m+1$.  Any sequence $w_m^n$ with $m > n$ represents the empty string
and is denoted by $\lambda$. The length of the empty string is
$\ell(\lambda) = 0$.

Given two finite sequences $w$ and $v$, we will denote by
$vw$ the sequence of length $\ell(v) + \ell(w) $ obtained by
concatenating the two strings.  In particular, $\lambda w = w\lambda=
w$.  The concatenation of sequences is also extended to the case in
which $v$ denotes a semi-infinite sequence, that is
$v=v_{-\infty}^{-1}$.

We say that the sequence $s$ is a \emph{suffix} of the sequence $w$ if
there exists a sequence $u$, with $\ell(u)\geq 1$, such that $w = us$.
In this case we write $s\prec w$. When $s\prec w$ or $s=w$ we write
$s\preceq w$. Given a sequence $w$ we denote by $\suf(w)$ the largest suffix of $w$.

In the sequel $A^j$ will denote the set of all sequences of length $j$
over $A$ and $A^*$ represents the set of all finite sequences, that is
\begin{equation*}
  A^* \,=\, \bigcup_{j=1}^{\infty}\,A^j.
\end{equation*}

\begin{definition} 
  A countable subset $\tau$ of $A^*$ is a \emph{tree} if no sequence
  $s \in \tau$ is a suffix of another sequence $w \in \tau$. This
  property is called the \emph{suffix property}.
\end{definition}

We define the \emph{height} of the tree $\tau$ as
\begin{equation*}
  h(\tau) = \sup\{\ell(w) : w\in\tau\}. 
\end{equation*}

In the case $h(\tau)<+\infty$ it follows that $\tau$ has a finite
number of sequences. In this case we say that $\tau$ is \emph{bounded}
and we will denote by $|\tau|$ the number of sequences in $\tau$. On
the other hand, if $h(\tau)=+\infty$ then $\tau$ has a countable
number of sequences. In this case we say that the tree $\tau$ is
\emph{unbounded}.
 
Given a tree $\tau$ and an integer $K$ we will denote by $\tau|_K$ the
tree $\tau$ \emph{truncated} to level $K$, that is
\begin{equation*}
  \tau|_K = \{w \in \tau\colon \ell(w) \le K\} \cup \{ w\colon \ell(w)=K 
  \text{ and } w\prec u, \text{ for some } u\in\tau\}.
\end{equation*}

We will say that a tree is \emph{irreducible} if no sequence can be
replaced by a suffix without violating the suffix property. This
notion was introduced in \cite{csiszar2006} and generalizes the
concept of complete tree. 
 
\begin{definition}\label{def:pct}
 A \emph{probabilistic context tree over $A$} is an
  ordered pair $(\tau,p)$ such that
\begin{enumerate}
\item $\tau$ is an irreducible tree;
\item $p = \{p(\cdot|w); w \in \tau\}$ is a family of transition
  probabilities over $A$.
\end{enumerate}
\end{definition}

Consider a stationary stochastic chain $(X_t)_{t\in\Z}$
over $A$. Given a sequence $w\in A^j$ we denote by 
\begin{equation*}
 p(w) \,=\,  \P(X_1^j = w)
\end{equation*}
the stationary probability of the cylinder defined by the sequence $w$.
If $p(w) > 0$ we write
\begin{equation*}
p(a|w) \,=\, \P ( X_0 =a \given X_{-j}^{-1}=w)\,.
\end{equation*}

\begin{definition}
  A sequence $w\in A^j$ is a \emph{context} for the process $(X_t)$ if
  $p(w)>0$ and for any semi-infinite sequence $x_{-\infty}^{-1}$ such
  that $w$ is a suffix of $x_{-\infty}^{-1}$ we have that
\begin{equation*}\label{eq:pt}
  \P ( X_0 =a \given X_{-\infty}^{-1}=x_{-\infty}^{-1}) \,=\, p(a|w),
  \quad\text{for all $a\in A$},
\end{equation*}
and no suffix of $w$ satisfies this equation.
\end{definition}

\begin{definition}
We say that the process $(X_t)$ is \emph{compatible} with the probabilistic context tree $(\tau,\bar p)$ if the following conditions are satisfied
\begin{enumerate}
\item $w\in\tau$ if and only if $w$ is a context for the process $(X_t)$. 
\item For any $w\in\tau$ and any $a\in A$, $\bar p(a|w) =  \P ( X_0 =a \given X_{-\ell(w)}^{-1}=w)$.
\end{enumerate}
\end{definition}

%

Define the sequence $(\alpha_k)_{k\in\N}$ as
\begin{align*}\label{alpha}
\alpha_0 &:=   \sum_{a\in A} \inf_{w\in\tau} \{\, p(a|w) \,\},\notag \\
\alpha_k &:= \inf_{u\in A^k}\;\sum_{a\in A}\;\inf_{w\in\tau, w\succ u}\{\,p(a|w)\,\}.
\end{align*}

From now on we will assume that the probabilistic context tree $(\tau,p)$ satisfies the following assumptions.
\begin{assumption}\label{ass1}
 Non-nullness, that is $\inf_{w\in\tau}\{p(a|w)\} > 0$ for any $a\in A$.  
 \end{assumption}

\begin{assumption}\label{ass2}
 Summability of the sequence $(1-\alpha_k), k\geq 0$. In this case denote by  
\[
\alpha := \sum_{k\in\N} (1-\alpha_k)\; < +\; \infty.
\]
\end{assumption}

For a probabilistic context tree satisfying Assumptions \ref{ass1} and \ref{ass2}, the maximal coupling argument used in \cite{fernandez2002}, or alternatively the perfect simulation scheme presented in \cite{comets2002}, imply the uniqueness of the law of the chain compatible with it.

Given an integer $k\geq 1$, we define
\begin{equation*}
\mathcal{C}_k = \{u\in\tau|_k\colon p(a|u)\neq p(a|\suf(u))\text{ for some }a\in A\}
\end{equation*}
and
\begin{equation*}
D_k = \min_{u\in \mathcal{C}_k}
\max_{a\in A} \{|p(a|u)-p(a|\suf(u))|\}.
\end{equation*}
We denote by  
\begin{equation*}\label{ek}
\epsilon_k = \min\{\,p(w)\colon \ell(w)\leq k \text{ and } p(w) > 0\,\}.
\end{equation*}

In what follows we will assume that $x_0, x_1, \dotsc, x_{n-1}$ is a
sample of the stationary stochastic chain $(X_t)$ compatible with the
probabilistic context tree $(\tau,p)$.

For any finite string $w$ with $\ell(w) \le n$, we denote by
$N_n(w)$ the number of occurrences of the string in the sample; that is
 \begin{equation*}
 \label{eq:Nn}
 N_n(w)=\sum_{t=0}^{n-\ell(w)}\Id\{X_t^{t+\ell(w)-1}=w\}.
 \end{equation*}

For any element $a \in A$ , the empirical transition probability
$\hat{p}_n(a|w)$ is defined by
\begin{equation}\label{phat}
\hat{p}_n(a|w)= \frac{N_n(wa)+1}{N_n(w\cdot)+|A|} \/.
\end{equation}
where 
\[
N_n(w\cdot)=\sum_{b \in A} N_n(wb)\,.
\]

This definition of $\hat{p}_n(a|w)$ is convenient because it is
asymptotically equivalent to $\frac{N_n(wa)}{N_n(w\cdot)}$ and it
avoids an extra definition in the case $N_n(w\cdot)=0$.

A variant of Rissanen's {\sl algorithm Context} is defined as follows.
First of all, let us define for any finite string $w\in A^*$:
\[ 
\Delta_n(w) = \max_{a\in A}
|\hat{p}_n(a|w)-\hat{p}_n(a|\suf(w))| \/.
\]
The $\Delta_n(w)$ operator computes a distance between the empirical transition probabilities 
associated to the sequence  $w$ and the one associated to the sequence $\suf(w)$.

\begin{definition}\label{estim-tree}
  Given $\delta > 0$ and $d < n$, the tree estimated with the algorithm Context is
  \begin{align*}
\hat{\tau}_n^{\delta,d} = \{w\in A_1^d\colon  &N_n(aw\cdot)>0,\Delta_n(a\,\suf(w)) >
\delta\text{ for some $a\in A$ and }\\
&\Delta_n(uw)\leq \delta \;\;\text{for all }u\in A_1^{d-\ell(w)}\text{ with }N_n(uw\cdot)\geq 1\},
\end{align*}
where $A_1^r$ denotes the set of all sequences of length at most $r$. In the case 
$\ell(w) = d$ we have $A_1^{d-\ell(w)} = \emptyset$.
\end{definition}

It is easy to see that $\hat{\tau}_n^{\delta,d}$ is an irreducible tree. Moreover, the way we defined 
$\hat{p}_n(\cdot|\cdot)$  in (\ref{phat})
associates a probability distribution to each sequence in $\hat{\tau}_n^{\delta,d}$.

The main result in this article is the following 

\begin{theorem}\label{expobounded}
 Let $(\tau,p)$ be a probabilistic context tree satisfying Assumptions \ref{ass1} and \ref{ass2}
 and let $(X_t)$ be a stationary stochastic chain compatible with $(\tau,p)$.
 Then for any integer $K$, any $d$ satisfying 
 \begin{equation}\label{d}
 d >  \max_{u\notin\tau\!\!,\,\ell(u)\leq K} \min\,\{k\colon \exists w\in\mathcal{C}_k, w\succ u\}
 \end{equation}
%
  any 
 $\delta < D_d$ and any
 \[
 n >  \frac{2(|A|+1)}{\min(\delta,D_d-\delta)\epsilon_d}+d
 \]
 we
  have that 
  \begin{equation*} 
  \P(\hat\tau_n^{\delta,d}|_K \neq \tau|_K)\;\leq\; 4\,\ee\, |A|^{d+2}\, \exp \bigl[- (n-d)\;
\frac{[\min(\frac{\delta}2,\frac{D_d-\delta}2)-\frac{|A|+1}{(n-d)\epsilon_d}]^2\epsilon_d^2C}{4|A|^2(d+1)},
\end{equation*}
  where
  \[
  C = \frac{\alpha_0}{8e(\alpha+\alpha_0)}.
  \]
\end{theorem}

As a consequence we obtain the following strong consistency result.

\begin{corollary}\label{main_cor}
 Under the conditions of Theorem~\ref{expobounded}
  we have 
  \begin{equation*}
  \hat\tau_n^{\delta,d}|_K = \tau|_K,
  \end{equation*}
eventually almost surely as $n\to+\infty$.
  \end{corollary}

\section{Exponential inequalities for empirical probabilities}\label{proof1} 

The main ingredient in the proof of Theorem~\ref{expobounded} is the following exponential upper bound

\begin{theorem}\label{estim1} 
  For any finite sequence $w$, any symbol $a\in A$ and any $t>0$ the
  following inequality holds
  \begin{equation*}\label{Nn2}
    \P(\,|N_n(wa)-(n-\ell(w))p(wa)|\,>\,t\,)\,\leq \,e^{\frac1e} 
    \exp \bigl[\frac{-t^2 C}{(n-\ell(w))\ell(wa)}\bigr],
  \end{equation*}
  where
  \begin{equation}\label{C}
    C = \frac{\alpha_0}{8e(\alpha+\alpha_0)}.
  \end{equation}
\end{theorem}

As a direct consequence of Theorem~\ref{estim1} we obtain the
following corollary.

\begin{corollary}\label{cor:estim1} 
  For any finite sequence $w$ with $p(w)>0$, any symbol $a\in A$,
  any $t>0$ and any $n > \frac{|A|+1}{tp(w)}+\ell(w)$ the following inequality holds
\begin{equation*}\label{expo2}
  \P\bigl(|\hat{p}_n(a|w)-p(a|w)|>t \bigl)\;\leq\;2\,|A|\,\ee \exp \bigl[- (n-\ell(w))\;\frac{[t-\frac{|A|+1}{(n-\ell(w))p(w)}]^2 p(w)^2C}{4|A|^2\ell(wa)}\bigl],
\end{equation*}
where $C$ is given by $(\ref{C})$.
\end{corollary}

To prove Theorem~\ref{estim1} we need a mixture property for processes compatible with a
probabilistic context tree $(\tau,p)$ satisfying Assumptions \ref{ass1} and \ref{ass2}. This 
is the content of the following lemma.

\begin{lemma}\label{mixing}
Let $(X_t)$ be a stationary stochastic chain compatible with the
  probabilistic context tree $(\tau,p)$ satisfying Assumptions \ref{ass1} and \ref{ass2}. Then, 
  there exists a summable sequence $\{\rho_l\}_{l\in\N}$, satisfying 
  \begin{equation}\label{rho}
  \sum_{l\in\N} \rho_l \;\leq\; 1+\frac{2 \alpha}{\alpha_0},
  \end{equation}  
   such that
  for any $i\geq 1$, any $k > i$, any $j\geq 1$ and any
  finite sequence $w_1^j$, the following inequality holds
  \begin{equation}\label{eqmixing}
    \sup_{x_1^{i}\in A^i} 
    |\P(X_{k}^{k+j-1}=w_1^j \given X_1^i=x_1^i)-p(w_1^j)|\;
    \leq \;  \sum_{l=0}^{j-1}\rho_{k-i-1+l}\,.
  \end{equation}
  \end{lemma}

\begin{proof} 
First note that  
  \begin{align*} \label{borne}
    \inf_{{u}\in A^\infty}\,
    \P(X_{k}^{k+j-1}=w_1^j \given X_{-\infty}^i&=u_{-\infty}^0x_1^i)  \,\leq \,
    \P(X_{k}^{k+j-1}=w_1^j \given X_1^i=x_1^i)\\[.2cm]
   & \;\leq \; \sup _{{u}\in A^\infty} 
    \,\P(X_{k}^{k+j-1}=w_1^j \given X_{-\infty}^i=u_{-\infty}^0x_1^i). 
  \end{align*} 
  where $A^\infty$ denotes the set of all semi-infinite sequences
  $u_{-\infty}^{0}$.  The reader can find a proof of the inequalities
  above in \cite[Proposition 3]{fernandez2002}. Using this fact and
  the condition of stationarity it is sufficient to prove that for any
  $k \geq 0$,
  \begin{equation*}
    \sup_{x\in A^\infty} 
    |\P(X_k^{k+j-1}=w_1^j \given X_{-\infty}^{-1}=x_{-\infty}^{-1})-p(w_1^j)|
    \leq \; \sum_{l=0}^{j-1}\rho_{k+l}.
  \end{equation*}
  Note that for all pasts $x_{-\infty}^{-1}$ we have 
  \begin{align*}
    \bigl| \P (X_k^{k+j-1}=w_1^j &\given
    X_{-\infty}^{-1}=x_{-\infty}^{-1})-p(w_1^j)\bigl|\\[.2cm]
    &= \; \Bigl|\int_{u\in A^\infty}
    \,\bigl[\P(X_k^{k+j-1}=w_1^j\given X_{-\infty}^{-1}=x_{-\infty}^{-1})\\
    &\mspace{140mu} - \P(X_k^{k+j-1}=w_1^j\given X_{-\infty}^{-1}=
    u_{-\infty}^{-1})\bigr] dp(u)\Bigl|\\
    &\leq\; \int_{u\in A^\infty}\,
    \bigl|\P(X_k^{k+j-1}=w_1^j\given X_{-\infty}^{-1}=x_{-\infty}^{-1})\\
    &\mspace{140mu}- \P(X_k^{k+j-1}=w_1^j\given
    X_{-\infty}^{-1}=u_{-\infty}^{-1}) \bigl| \,dp(u).
 \end{align*}
Therefore, applying 
the loss of memory property proved in \cite[Corollary 4.1]{comets2002} we have that
  \begin{equation*}\label{eqmixing2}
  \bigl|\P(X_k^{k+j-1}=w_1^j\given X_{-\infty}^{-1}=x_{-\infty}^{-1})
   - \P(X_k^{k+j-1}=w_1^j\given 
    X_{-\infty}^{-1}=u_{-\infty}^{-1}) \bigl| \; \leq \;   \sum_ {l=0}^{j-1} \rho_{k+l},
 \end{equation*}
 where $\rho_m$ is defined as the probability of return to the origin at time $m$ of the 
   Markov chain on $\N$ starting at time zero at the origin and having transition probabilities
  \begin{align}\label{hoc}
  p(x,y) = \begin{cases}
  \alpha_x, & \text{ if y = x+1},\\  
  1-\alpha_x, & \text{ if y=0},\\
  0, & \text{ otherwise}.
  \end{cases}
  \end{align}
This concludes the proof of (\ref{eqmixing}). To prove 
(\ref{rho}), let $(Z_n)$ be the Markov chain with probability transitions given by 
(\ref{hoc}). By definition we have 
 \begin{align*}\label{rho-alpha}
    \prod_{l\geq 1}(1-\rho_l)  \; & = \prod_{l\geq 1}\,\sum_{j=1}^l \P(Z_{l}=j\given Z_{l-1}=j-1)\P(Z_{l-1}=j-1)\\
    & \geq \;  \prod_{l\geq 1}  \alpha_{l-1} \prod_{i=0}^{l-2}\alpha_{i}\,\geq  \prod_{l\geq 0}  \alpha_{l}^2. 
\end{align*}
From this, using the inequality $x \,\leq\, -\ln (1-x)\, \leq \,\frac{x}{1-c}$ which holds for any $x\in(-1,c\,]$, it follows  that
  \[
  \sum_{l\geq 1}\rho_l\;\leq\; -2 \sum_{l\geq 0} \log \alpha_l \;\leq \;2 \sum_{l\geq 0} \frac{1-\alpha_l}{\alpha_0}.
  \]
This concludes the proof of the lemma.
\end{proof}

We are now ready to prove Theorem~\ref{estim1}. 


\begin{proof}[Proof of Theorem~\ref{estim1}]
  Let $w$ be a finite sequence and $a$ any symbol in $A$. 
  Define the random variables
  \begin{equation*}
  U_j = \Id\{X_j^{j+\ell(w)} = wa\} - p(wa),
  \end{equation*}
  for $j=0,\dotsc,n-\ell(wa)$. Then, using \cite[Proposition~4]{DD} 
  we have that, for any $p\geq 2$
  \begin{align*}
    \Vert N_n&(wa)- (n-\ell(w))p(wa)\Vert_p \\[.1cm]
      &\,\leq\, \Bigl( 2p
    \sum_{i=0}^{n-\ell(wa)}\sum_{k=i}^{n-\ell(wa)} 
    \Vert \,\E(U_k\given U_0,\dotsc,U_i)\,\Vert_{\infty}\Bigr)^{\frac12}\\
    &\,\leq\, \Bigl( 2p
    \sum_{i=0}^{n-\ell(wa)}\sum_{k=i}^{n-\ell(wa)} \sup_{u\in A^{i+\ell(wa)}}
    |\P(X_k^{k+\ell(w)}=wa\given
    X_0^{i+\ell(w)}= u)-p(wa)|\Bigr)^{\frac12}\\
    &\,\leq \,\Bigl( 2p\,\ell(wa)(n-\ell(w))\frac{2(\alpha+\alpha_0)}{\alpha_0}\Bigr)^{\frac12}.
  \end{align*}
  Then, as in \cite[Proposition~5]{DP} we also obtain that, for any $t>0$,
  \begin{equation*}
    \P(|N_n(wa)-(n-\ell(w))\/p(wa)|>t)\, \leq \,e^{\frac1e}\, 
    \exp \bigl[\frac{-t^2C}{(n-\ell(w))\ell(wa)}\bigl]\/,
  \end{equation*}
  where
  \[
  C = \frac{\alpha_0}{8e(\alpha+\alpha_0)}.
  \]
\end{proof}

\begin{proof}[Proof of Corollary~\ref{cor:estim1}]
  First observe that
  \begin{equation*}
  \Bigl|\, p(a|w) - \frac{(n-\ell(w))p(wa)+1}{(n-\ell(w)) p(w)+|A|}\, \Bigr|\;\leq\;
  \frac{|A|+1}{(n-\ell(w))p(w)}\/.
\end{equation*}
  Then, for all $n\geq (|A|+1)/tp(w) + \ell(w)$ we have that
 \begin{align*}
    \P\bigl(\,\bigl|\hat{p}_n(a|w) &- p(a|w)\bigr|>t\,\bigl)\\
    &\leq \;\P\bigl(\,\Bigl|\frac{N_n(wa)+1}{N_n(w\cdot)+|A|} - \frac{(n-\ell(w))p(wa)+1}{(n-\ell(w))p(w)+|A|} \Bigr|> t - \frac{|A|+1}{(n-\ell(w))p(w)}\, \bigr)\\
  \end{align*}
Denote by $t'=  t - (|A|+1)/(n-\ell(w))p(w)$. Then 
  \begin{align*}
  \P\bigl(\,\Bigl|\frac{N_n(wa)+1}{N_n(w\cdot)+|A|}& - \frac{(n-\ell(w))p(wa)+1}{(n-\ell(w))
  p(w)+|A|}\Bigr|> t'\,
  \bigr)\\[.2cm]
   &\mspace{-60mu}\leq \;\P\bigl(\bigl|N_n(wa)-(n-\ell(w))p(wa)\bigr|>
    \frac{t'}2 [(n-\ell(w))p(w)+|A|]  \bigl) \\[.2cm]
  &\mspace{-60mu} \quad + \sum_{b\in A}\;\P\bigl(\bigl|N_n(wb)-(n-\ell(w))p(wb)\bigr|>
\frac{t'}{2|A|} [ (n-\ell(w))p(w)+|A|] \bigl).
\end{align*}
Now, we can apply Theorem~\ref{estim1} to bound above the last sum by
\begin{equation*}
2\,|A|\, \ee\, \exp \bigl[- (n-\ell(w))\;\frac{[t-\frac{|A|+1}{(n-\ell(w))p(w)}]^2 p(w)^2C}{4|A|^2\ell(wa)}\bigl],
\end{equation*}
 where
 \[
 C = \frac{\alpha_0}{8e(\alpha+\alpha_0)}.
 \]
 This finishes the proof of the corollary.
  \end{proof}

\section{Proof of the main results}\label{proofmain}

\begin{proof}[Proof of Theorem~\ref{expobounded}]

\noindent Define
\[
O_{n}^{\delta,d}=\bigcup_{\substack{w \in \tau\\[.1cm]\ell(w)< K}}\bigcup_{uw\in\hat\tau_n^{\delta,d}} \{ \Delta_n(uw) >\delta \}
\/,
\]
and
\[
U_{n}^{\delta,d}=\bigcup_{\substack{w \in \hat\tau_n^{\delta,d}\\[.1cm]\ell(w) < K}}\bigcap_{uw\in\tau|_d} \{\Delta_n(uw) \leq \delta\}.
\]
Then, if  $d < n$ we have that
\[
\{\hat{\tau}_n^{\delta,d}|_K \neq \tau|_K\} = O_{n}^{\delta,d}\cup U_{n}^{\delta,d}.
\]
The result follows from a succession of lemmas.


\begin{lemma}\label{lemOn}
For any $n > \frac{2(|A|+1)}{\delta\epsilon_d}+d$, for any $w \in \tau$ with $\ell(w)<K$ and for any $uw\in\hat{\tau}_n^{\delta,d}$ we have that
  \begin{align*}
  \P( \Delta_n(uw) >\delta  )\;\leq \;  4\, |A|^2\,\ee \exp \bigl[- (n-d)\;\frac{[\frac{\delta}2-\frac{|A|+1}{(n-d)\epsilon_d}]^2 \epsilon_d^2C}{4|A|^2(d+1)}\bigl],
  \end{align*}
  where  $C$ is given by (\ref{C}).
\end{lemma}

\begin{proof}
 Recall that 
  \[
  \Delta_n(uw) = \max_{a\in A} |\hat{p}_n(a|uw)-\hat{p}_n(a|\suf(uw))|.
  \] 
  Note that the fact $w\in\tau$ implies that
  for any finite sequence $u$ with $p(u)>0$ and any symbol
   $a\in A$ we have $p(a|w)=p(a|uw)$.
  Hence,
  \begin{align*}
    \P(\Delta_n(uw) >\delta)\;\leq\; \sum_{a\in A}\,\bigl[&\, 
    \P\bigl( |\hat{p}_n(a|w)-p(a|w)|>\frac\delta{2}\bigr)\\
    & +\P\bigl( |\hat{p}_n(a|uw)-p(a|uw)|>\frac\delta{2}\bigr)\bigr].
  \end{align*}
  Using Corollary~\ref{cor:estim1} we can bound above the right hand side of the last inequality by 
  \[
4\, |A|^2\,\ee \exp \bigl[- (n-d)\;\frac{[\frac{\delta}2-\frac{|A|+1}{(n-d)\epsilon_d}]^2 \epsilon_d^2C}{4|A|^2(d+1)}\bigl],
   \]
   where $C$ is given by (\ref{C}).
\end{proof}


\begin{lemma}\label{lemUn}
  For any $n > \frac{2(|A|+1)}{(D_d -\delta)\epsilon_d}+d$ and for any  
   $w\in\hat{\tau}_n^{\delta,d}$ with $\ell(w)<K$ we have that
  \[
  \P(\bigcap_{uw\in\tau|_d} \{\Delta_n(uw) \leq \delta\})\,\leq 
  4\,|A|\, \ee \exp \bigl[- (n-d)\;\frac{[\frac{D_d-\delta}{2}-\frac{|A|+1}{(n-d)\epsilon_d}]^2 \epsilon_d^2C}{4|A|^2(d+1)}\bigl],
    \]
    where $C$ is given by $(\ref{C})$.
\end{lemma}
\begin{proof}
As $d$ satisfies (\ref{d}) there exists $\bar{uw}\in\tau|_d$ such that
  $p(a|\bar{uw})\neq p(a|\suf(\bar{uw}))$ for some $a\in A$. Then 
  \begin{equation*}
   \P(\bigcap_{uw\in\tau|_d} \{\Delta_n(uw) \leq \delta\})\,\leq\,  \P(\Delta_n(\bar{uw}) \leq \delta).
  \end{equation*}
Observe that for  any $a\in A$,
  \begin{align*}
    |\hat{p}_n(a|\suf(\bar{uw}))-\hat{p}_n(a|&\bar{uw})|\;\geq \; |p(a|\suf(\bar{uw}))-p(a|\bar{uw})| \\
  &  - |\hat{p}_n(a|\suf(\bar{uw}))-p(a|\suf(\bar{uw}))| 
   - |\hat{p}_n(a|\bar{uw})-p(a|\bar{uw})|.
  \end{align*}
   Hence, we have that for any $a\in A$ 
  \[
  \Delta_n(\bar{uw}) \, \geq\, D_d - |\hat{p}_n(a|\suf(\bar{uw}))-p(a|\suf(\bar{uw}))| - |\hat{p}_n(a|\bar{uw})-p(a|\bar{uw})|\/.
  \]
  Therefore,
  \begin{eqnarray*}
    \P(\Delta_n(\bar{uw}) \leq \delta)&\leq &\P\bigl(\,\bigcap_{a\in A}
    \{\,|\hat{p}_n(a|\suf(\bar{uw}))-p(a|\suf(\bar{uw}))|\geq\frac{D_d-\delta}2\,\}\,\bigl) \\
    & & + \P\bigl(\,\bigcap_{a\in A}
    \{\,|\hat{p}_n(a|\bar{uw})-p(a|\bar{uw})|\geq\frac{D_d-\delta}2\,\}\,\bigl)\/.
  \end{eqnarray*}
  As $\delta < D_d$ and $n > \frac{2(|A|+1)}{(D_d -\delta)\epsilon_d}+d$ we can use Corollary~\ref{cor:estim1} to bound above the right hand side of this 
  inequality by
  \[
  4\,|A|\, \ee \exp \bigl[- (n-d)\;\frac{[\frac{D_d-\delta}{2}-\frac{|A|+1}{(n-d)\epsilon_d}]^2\epsilon_d^2C}{4|A|^2(d+1)}\bigl],
\]
where $C$ is given by $(\ref{C})$.
  This concludes the proof of the lemma.
\end{proof}

\noindent  Now we can finish the proof of Theorem~\ref{expobounded}. We have that
\[
\P(\hat{\tau}_n^{\delta,d}|_K \neq \tau|_K) = \P(O_{n}^{\delta,d})+\P(U_{n}^{\delta,d}).
\]
Using the definition of $O_{n}^{\delta,d}$ and $U_{n}^{\delta,d}$ we have that 
\[
\P(\hat{\tau}_n^{\delta,d}|_K \neq \tau|_K)\leq\! \sum_{\substack{w \in
    \tau\\ \ell(w)<K}}\sum_{uw\in\hat{\tau}_n^{\delta,d}} \P(  \Delta_n(uw) > \delta) + \!\sum_{\substack{w \in
    \hat\tau_n^{\delta,d}\\\ell(w)<K}}\P ( \bigcap_{uw\in\tau|_d}  \Delta_n(uw) \leq \delta ) .
\]
Applying Lemma~\ref{lemOn} and Lemma~\ref{lemUn} we can bound above the last expression by  
\begin{equation*}
\P(\hat{\tau}_n^{\delta,d}|_K \neq \tau|_K)\,\leq\, 4\,\ee\, |A|^{d+2}\, \exp \bigl[- (n-d)\;
\frac{[\min(\frac{\delta}2,\frac{D_d-\delta}2)-\frac{|A|+1}{(n-d)\epsilon_d}]^2\epsilon_d^2C}{4|A|^2(d+1)}\bigl],
\end{equation*}
where $C$ is given by (\ref{C}).
We conclude the proof of Theorem~\ref{expobounded}. 
\end{proof}

\begin{proof}[Proof of Corollary~\ref{main_cor}]
It follows from Theorem~\ref{expobounded}, using 
the first Borel-Cantelli Lem\-ma and
the fact that the bounds for the error estimation of the context tree are summable in $n$ for a fixed $d$ satisfying (\ref{d}) and $\delta<D_d$. 
\end{proof}

\section{Final remarks}

The present paper presents an upper bound for the rate of convergence 
of a version of the algorithm Context, for unbounded context trees.
This generalizes previous results obtained in \cite{galves2006} for the case of bounded variable memory 
processes. We obtain an exponential bound for the probability of incorrect estimation of the truncated context tree, when the estimator is given by Definition~(\ref{estim-tree}). 
Note that the definition of the context tree estimator depends on the parameter $\delta$, and this parameter appears in the exponent of the upper bound.  
To assure the consistency of the estimator we need to choose a $\delta$ sufficiently small, depending on the
transition probabilities of the process. Therefore, our estimator is not universal, in the sense that for any fixed $\delta$ it fails to be consistent for any process having $D_d < \delta$. The same happens with the parameter $d$. In order to choose $\delta$ and $d$ not depending on the process, we can allow these parameters to be a function of $n$, in such a way $\delta_n$ goes to zero and $d_n$ goes to $+\infty$ as $n$ diverges. When we do this, we loose the exponential property of the upper bound. 

As an anonymous referee has pointed out, Finesso et al. \cite{finesso1996} proved that in the simpler case of estimating the order of a Markov chain, it is not possible to obtain pure exponential bounds for the overestimation event with a universal estimator. The above discussion illustrates this fact.


\section{Acknowledgments}
We thank Pierre Collet, Imre Csisz\'{a}r, Nancy Garcia, Aur\'elien Garivier, Bezza Hafidi, V\'{e}ronique Maume-Deschamps, Eric Moulines, Jorma Rissanen and Bernard Schmitt for many discussions on the subject. We also thank an anonymous referee that attracted our attention to the interesting paper \cite{finesso1996}. 

\bibliography{./referencias}  
\bibliographystyle{plain}

\end{document}